\input amstex
\documentstyle{amsppt}
%----------------------------------------------------------------
% Title:     A fast modulo primes algorithm for searching perfect 
%            cuboids and its implementation.
% Authors:   R.R.Gallyamov, I.R.Kadyrov, D.D.Kashelevskiy,
%            N.G.Kutlugallyamov, R.A.Sharipov
% Comments:  AmSTeX, 11 pages, amsppt style
% MSC-class: 11D41, 11D72, 68U99, 65-04
%----------------------------------------------------------------
%           Replacement for output macro definition
%
\catcode`@=11
\redefine\output@{%
  \def\break{\penalty-\@M}\let\par\endgraf
  \ifodd\pageno\global\hoffset=105pt\else\global\hoffset=8pt\fi  
  \shipout\vbox{%
    \ifplain@
      \let\makeheadline\relax \let\makefootline\relax
    \else
      \iffirstpage@ \global\firstpage@false
        \let\rightheadline\frheadline
        \let\leftheadline\flheadline
      \else
        \ifrunheads@ %\let\makefootline\relax
        \else \let\makeheadline\relax
        \fi
      \fi
    \fi
    \makeheadline \pagebody \makefootline}%
  \advancepageno \ifnum\outputpenalty>-\@MM\else\dosupereject\fi
}
\def\Beta{\mathchar"0\hexnumber@\rmfam 42}
\catcode`\@=\active
%----------------------------------------------------------------
\nopagenumbers
\chardef\textvolna='176

\chardef\bigalpha='013
\def\negskp{\hskip -2pt}

\chardef\degree="5E
\def\compos{\,\raise 1pt\hbox{$\sssize\circ$} \,}

\def\vtrule{\vrule height 12pt depth 6pt}
%\font\eightrm=cmr8
%\def\LT{\operatorname{\text{\eightrm LT}}}
%\def\LM{\operatorname{\text{\eightrm LM}}}
%\def\LC{\operatorname{\text{\eightrm LC}}}
%\accentedsymbol\hatgamma{\kern 2pt\hat{\kern -2pt\gamma}}
%\accentedsymbol\checkgamma{\kern 2.5pt\check{\kern -2.5pt\gamma}}
\def\blue#1{#1}
\def\red#1{#1}
\gdef\darkred#1{#1}
\gdef\green#1{#1}
\catcode`#=11\def\diez{#}\catcode`#=6
\catcode`&=11\catcode`&=4
\catcode`_=11\def\podcherkivanie{_}\catcode`_=8
\catcode`\^=11\catcode`\^=7
\catcode`~=11\def\volna{~}\catcode`~=\active
\def\mycite#1{\cite{\blue{#1}}\immediate\special{ps:
     ShrHPSdict begin /ShrBORDERthickness 0 def}}
\def\myciterange#1#2#3#4{\cite{\blue{#2#3#4}}\immediate\special{ps:
     ShrHPSdict begin /ShrBORDERthickness 0 def}}
\def\mytag#1{%
    \tag#1}
\def\mythetag#1{\thetag{\blue{#1}}\immediate\special{ps:
     ShrHPSdict begin /ShrBORDERthickness 0 def}}
\def\myrefno#1{\no#1}
\def\myhref#1#2{\blue{#2}\immediate\special{ps:
     ShrHPSdict begin /ShrBORDERthickness 0 def}}
\def\myEarXivlink{\myhref{http://arXiv.org}{http:/\negskp/arXiv.org}}

\def\mytheorem#1{\csname proclaim\endcsname{Theorem #1}}
\def\mytheoremwithtitle#1#2{\csname proclaim\endcsname{Theorem #1#2}}
\def\mythetheorem#1{\blue{#1}\immediate\special{ps:
     ShrHPSdict begin /ShrBORDERthickness 0 def}}
\def\mylemma#1{\csname proclaim\endcsname{Lemma #1}}
\def\mylemmawithtitle#1#2{\csname proclaim\endcsname{Lemma #1#2}}

\def\mycorollary#1{\csname proclaim\endcsname{Corollary #1}}

\def\mydefinition#1{\definition{Definition #1}}

\def\myconjecture#1{\csname proclaim\endcsname{Conjecture #1}}
\def\myconjecturewithtitle#1#2{\csname proclaim\endcsname{Conjecture #1#2}}
\def\mytheconjecture#1{\blue{#1}\immediate\special{ps:
     ShrHPSdict begin /ShrBORDERthickness 0 def}}
\def\myproblem#1{\csname proclaim\endcsname{Problem #1}}
\def\myproblemwithtitle#1#2{\csname proclaim\endcsname{Problem #1#2}}

%----------------------------------------------------------------
% Cyrillic fonts definition

%\font\tencyr=wncyr10
%----------------------------------------------------------------
\pagewidth{360pt}
\pageheight{606pt}
\topmatter
\title
A fast modulo primes algorithm for searching perfect cuboids
and its implementation.
\endtitle
\rightheadtext{A fast modulo primes algorithm \dots}
\author
R\.\,R\.\,Gallyamov, I\.\,R\.\,Kadyrov, D\.\,D\.\,Kashelevskiy,\\
N\.\,G\.\,Kutlugallyamov, R\.\,A\.\,Sharipov
\endauthor
\leftheadtext{R\.\,A\.\,Sharipov and students}
\address Bashkir State University, 32 Zaki Validi street, 450074 Ufa, Russia
\endaddress
\email
\myhref{mailto:r-sharipov\@mail.ru}{r-sharipov\@mail.ru}
\endemail
\abstract
     A perfect cuboid is a rectangular parallelepiped whose all linear extents
are given by integer numbers, i\.\,e\. its edges, its face diagonals, and its
space diagonal are of integer lengths. None of perfect cuboids is known thus 
far. Their non-existence is also not proved. This is an old unsolved mathematical 
problem.\par 
    Three mathematical propositions have been recently associated with the cuboid 
problem. They are known as three cuboid conjectures. These three conjectures specify 
three special subcases in the search for perfect cuboids. The case of the second
conjecture is associated with solutions of a tenth degree Diophantine equation. In 
the present paper a fast algorithm for searching solutions of this Diophantine 
equation using modulo primes seive is suggested and its implementation on 32-bit 
Windows platform with Intel-compatible processors is presented. 
\endabstract
\subjclassyear{2000}
\subjclass 11D41, 11D72, 68U99, 65-04\endsubjclass
\endtopmatter
%\loadbold
%\loadeufb
\TagsOnRight
\document
%\input countstyle

%\special{header=resource.eps}
\head
1. Introduction.
\endhead
\myconjecturewithtitle{1.1}{\rm\ (Second cuboid conjecture)} For any two positive coprime 
integer numbers $p\neq q$ the tenth-degree polynomial
$$
\gathered
Q_{pq}(t)=t^{10}+(2\,q^{\kern 0.7pt 2}+p^{\kern 1pt 2})\,(3\,q^{\kern 0.7pt 2}-2\,p^{\kern 1pt 2})\,t^8
+(q^{\kern 0.5pt 8}+10\,p^{\kern 1pt 2}\,q^{\kern 0.5pt 6}+\\
+\,4\,p^{\kern 1pt 4}\,q^4-14\,p^{\kern 1pt 6}\,q^{\kern 0.7pt 2}+p^{\kern 1pt 8})\,t^6
-p^{\kern 1pt 2}\,q^{\kern 0.7pt 2}\,(q^{\kern 0.5pt 8}-14\,p^{\kern 1pt 2}\,q^{\kern 0.5pt 6}+4\,p^{\kern 1pt 4}\,q^4+\\
+\,10\,p^{\kern 1pt 6}\,q^{\kern 0.7pt 2}+p^{\kern 1pt 8})\,t^4
-p^{\kern 1pt 6}\,q^{\kern 0.5pt 6}\,(q^{\kern 0.7pt 2}
+2\,p^{\kern 1pt 2})\,(3\,p^{\kern 1pt 2}-2\,q^{\kern 0.7pt 2})\,t^2
-q^{\kern 0.7pt 10}\,p^{\kern 1pt 10}
\endgathered\quad
\mytag{1.1}
$$
is irreducible over the ring of integers $\Bbb Z$.
\endproclaim
\mytheorem{1.1} A perfect cuboid associated with the polynomial  \mythetag{1.1} does exist 
if and only if for some positive coprime integer numbers $p\neq q$ the Diophantine equation $Q_{pq}(t)=0$ has a positive solution $t$ obeying the inequalities
$$
\xalignat 4
&t>p^{\kern 1pt 2},
&&t>p\,q,
&&t>q^{\kern 0.7pt 2},
&&(p^{\kern 1pt 2}+t)\,(p\,q+t)>2\,t^2.
\quad
\endxalignat
$$
\endproclaim
    Theorem~\mythetheorem{1.1} can be found in \mycite{1}. It stems from the 
results of \mycite{2} and \mycite{3}. As for the perfect cuboid problem itself,
it has a long history reflected in \myciterange{4}{4}{--}{51}. There are also 
two series of ArXiv publications. The first of them \myciterange{52}{52}{--}{54}
continues the research on cuboid conjectures. The second one 
\myciterange{55}{55}{--}{67} relates perfect cuboids with multisymmetric 
\pagebreak polynomials.\par
     The scope of perfect cuboids in the case of the second cuboid conjecture 
is restricted by the following theorem derived from \mycite{1}.\par
\mytheorem{1.2} In the case of the second cuboid conjecture there are no perfect
cuboids outside the region given by the inequalities 
$$
\hskip -2em
\min\left(\root{\raise 1pt\hbox{$\ssize 3\kern -1pt $}}\of{\frac{p}{9}\,},
\,\frac{p}{59}\right)\leqslant q\leqslant{59\,p}.
\mytag{1.2}
$$
\endproclaim
In \mycite{1} the region given by the inequalities \mythetag{1.2} was presented 
as a union of two regions which were called the linear and the nonlinear regions 
respectively. In this paper we present an algorithm for searching cuboids in the
region \mythetag{1.2}. 
\head
2. A modulo primes seive.
\endhead
     Let $p$, $q$, and $t$ be a triple of integer numbers satisfying the Diophantine
equation $Q_{pq}(t)=0$ with the polynomial \mythetag{1.1} and let $r$ be some prime
number. Then we can pass from $\Bbb Z$ to the quotient ring $\Bbb Z_r=\Bbb Z/r\Bbb Z$
and denote
$$
\xalignat 3
&\hskip -2em
\tilde p=p\!\!\mod r,
&&\tilde q=q\!\!\mod r,
&&\tilde t=t\!\!\mod r.
\mytag{2.1}
\endxalignat
$$
The numbers $\tilde p$, $\tilde q$, and $\tilde t$ are interpreted as division remainders
after dividing $p$, $q$, and $t$ by the prime number $r$. They obey the quotient equation
$$
\hskip -2em
Q_{\tilde p\tilde q}(\tilde t)\!\!\!\!\mod r=0. 
\mytag{2.2}
$$
Once $r$ is given there are only a finite number of remainders \mythetag{2.1}:
$$
\xalignat 3
&\tilde p=0,\,\ldots,\,r-1,
&&\tilde q=0,\,\ldots,\,r-1,
&&\tilde t=0,\,\ldots,\,r-1.
\quad
\endxalignat
$$
The values in the left hand side of the equation \mythetag{2.2} for them can be precomputed.
They can be either zero or nonzero modulo $r$. We can use them as a fast computed test for
sweeping away those values of $p$, $q$, and $t$, where $Q_{pq}(t)\neq 0$.
\mydefinition{1.1} A pair of integer numbers $0\leqslant\tilde p\leqslant r-1$ and
$0\leqslant\tilde q\leqslant r-1$ is called solvable modulo $r$ if there is at least one
integer number $0\leqslant\tilde t\leqslant r-1$ such that $Q_{\tilde p\tilde q}(\tilde t)\!\!\!\mod r=0$. Otherwise it is called unsolvable. 
\enddefinition 
     We can represent solvable and unsolvable pairs in the form of bit-arrays $u_r$:
$$
\hskip -2em
u_r(\tilde p,\tilde q)
=\cases 0 &\text{if $(\tilde p,\tilde q)$ is solvable;}\\
1 &\text{if $(\tilde p,\tilde q)$ is unsolvable.}
\endcases
\mytag{2.3}
$$
The value $u_r(\tilde p,\tilde q)$ of the function \mythetag{2.3} is called the
unsolvability bit. Bit-arrays of the form \mythetag{2.3} can be stored as tables.
For $r=2$ this table looks like 
$$
\def\vtrule{\vrule height 12pt depth 6pt}
\hskip -2em
\vcenter{\vbox{\hsize 4cm
\offinterlineskip\settabs\+\indent
\vtrule
\hskip 1.1cm &\vtrule
\hskip 1.1cm &\vtrule
\hskip 1.1cm &\vtrule
\cr\hrule 
\+\vtrule
\hss\,$u_2(\tilde p,\tilde q)$ &\vtrule
\hss\ $\tilde q=0$\hss &\vtrule
\hss\ $\tilde q=1$\hss &\vtrule\cr\hrule
\+\vtrule
\hss\ p=0\hss&\vtrule
\hss\ \ 0\hss&\vtrule
\hss\ \ 0\hss&\vtrule\cr\hrule
\+\vtrule
\hss\ p=1\hss&\vtrule
\hss\ \ 0\hss&\vtrule
\hss\ \ 0\hss&\vtrule\cr\hrule
}}
\mytag{2.4}
$$
As we see in \mythetag{2.4}, the values of the function $u_2(\tilde p,\tilde q)$
are identically zero. The same is true for the functions $u_3(\tilde p,\tilde q)$
$u_5(\tilde p,\tilde q)$, and $u_7(\tilde p,\tilde q)$ associated with the prime
numbers $r=3$, $r=5$, and $r=7$. The case of $r=11$ is different: 
$$
\hskip -2em
\vcenter{\vbox{\hsize 12cm
\offinterlineskip\settabs\+\indent
\vtrule
\hskip 0.7cm &\vtrule
\hskip 0.7cm &\vtrule
\hskip 0.7cm &\vtrule
\hskip 0.7cm &\vtrule
\hskip 0.7cm &\vtrule
\hskip 0.7cm &\vtrule
\hskip 0.7cm &\vtrule
\hskip 0.7cm &\vtrule
\hskip 0.7cm &\vtrule
\hskip 0.7cm &\vtrule
\hskip 0.7cm &\vtrule
\hskip 0.9cm &\vtrule
\cr\hrule 
\+\vtrule
\hss\,\ $u_{11}$\hss&\vtrule
\hss\ \ 0\hss&\vtrule
\hss\ \ 1\hss&\vtrule
\hss\ \ 2\hss&\vtrule
\hss\ \ 3\hss&\vtrule
\hss\ \ 4\hss&\vtrule
\hss\ \ 5\hss&\vtrule
\hss\ \ 6\hss&\vtrule
\hss\ \ 7\hss&\vtrule
\hss\ \ 8\hss&\vtrule
\hss\ \ 9\hss&\vtrule
\hss\ \ 10\hss&\vtrule
\cr\hrule
\+\vtrule
\hss\ p=0\hss&\vtrule
\hss\ \ 0\hss&\vtrule
\hss\ \ 0\hss&\vtrule
\hss\ \ 0\hss&\vtrule
\hss\ \ 0\hss&\vtrule
\hss\ \ 0\hss&\vtrule
\hss\ \ 0\hss&\vtrule
\hss\ \ 0\hss&\vtrule
\hss\ \ 0\hss&\vtrule
\hss\ \ 0\hss&\vtrule
\hss\ \ 0\hss&\vtrule
\hss\ \ 0\hss&\vtrule
\cr\hrule

\+\vtrule
\hss\ p=1\hss&\vtrule
\hss\ \ 0\hss&\vtrule
\hss\ \ 0\hss&\vtrule
\hss\ \ 1\hss&\vtrule
\hss\ \ 1\hss&\vtrule
\hss\ \ 1\hss&\vtrule
\hss\ \ 1\hss&\vtrule
\hss\ \ 1\hss&\vtrule
\hss\ \ 1\hss&\vtrule
\hss\ \ 1\hss&\vtrule
\hss\ \ 1\hss&\vtrule
\hss\ \ 0\hss&\vtrule
\cr\hrule
\+\vtrule
\hss\ p=2\hss&\vtrule
\hss\ \ 0\hss&\vtrule
\hss\ \ 1\hss&\vtrule
\hss\ \ 0\hss&\vtrule
\hss\ \ 1\hss&\vtrule
\hss\ \ 1\hss&\vtrule
\hss\ \ 1\hss&\vtrule
\hss\ \ 1\hss&\vtrule
\hss\ \ 1\hss&\vtrule
\hss\ \ 1\hss&\vtrule
\hss\ \ 0\hss&\vtrule
\hss\ \ 1\hss&\vtrule
\cr\hrule
\+\vtrule
\hss\ p=3\hss&\vtrule
\hss\ \ 0\hss&\vtrule
\hss\ \ 1\hss&\vtrule
\hss\ \ 1\hss&\vtrule
\hss\ \ 0\hss&\vtrule
\hss\ \ 1\hss&\vtrule
\hss\ \ 1\hss&\vtrule
\hss\ \ 1\hss&\vtrule
\hss\ \ 1\hss&\vtrule
\hss\ \ 0\hss&\vtrule
\hss\ \ 1\hss&\vtrule
\hss\ \ 1\hss&\vtrule
\cr\hrule
\+\vtrule
\hss\ p=4\hss&\vtrule
\hss\ \ 0\hss&\vtrule
\hss\ \ 1\hss&\vtrule
\hss\ \ 1\hss&\vtrule
\hss\ \ 1\hss&\vtrule
\hss\ \ 0\hss&\vtrule
\hss\ \ 1\hss&\vtrule
\hss\ \ 1\hss&\vtrule
\hss\ \ 0\hss&\vtrule
\hss\ \ 1\hss&\vtrule
\hss\ \ 1\hss&\vtrule
\hss\ \ 1\hss&\vtrule
\cr\hrule
\+\vtrule
\hss\ p=5\hss&\vtrule
\hss\ \ 0\hss&\vtrule
\hss\ \ 1\hss&\vtrule
\hss\ \ 1\hss&\vtrule
\hss\ \ 1\hss&\vtrule
\hss\ \ 1\hss&\vtrule
\hss\ \ 0\hss&\vtrule
\hss\ \ 0\hss&\vtrule
\hss\ \ 1\hss&\vtrule
\hss\ \ 1\hss&\vtrule
\hss\ \ 1\hss&\vtrule
\hss\ \ 1\hss&\vtrule
\cr\hrule
\+\vtrule
\hss\ p=6\hss&\vtrule
\hss\ \ 0\hss&\vtrule
\hss\ \ 1\hss&\vtrule
\hss\ \ 1\hss&\vtrule
\hss\ \ 1\hss&\vtrule
\hss\ \ 1\hss&\vtrule
\hss\ \ 0\hss&\vtrule
\hss\ \ 0\hss&\vtrule
\hss\ \ 1\hss&\vtrule
\hss\ \ 1\hss&\vtrule
\hss\ \ 1\hss&\vtrule
\hss\ \ 1\hss&\vtrule
\cr\hrule
\+\vtrule
\hss\ p=7\hss&\vtrule
\hss\ \ 0\hss&\vtrule
\hss\ \ 1\hss&\vtrule
\hss\ \ 1\hss&\vtrule
\hss\ \ 1\hss&\vtrule
\hss\ \ 0\hss&\vtrule
\hss\ \ 1\hss&\vtrule
\hss\ \ 1\hss&\vtrule
\hss\ \ 0\hss&\vtrule
\hss\ \ 1\hss&\vtrule
\hss\ \ 1\hss&\vtrule
\hss\ \ 1\hss&\vtrule
\cr\hrule
\+\vtrule
\hss\ p=8\hss&\vtrule
\hss\ \ 0\hss&\vtrule
\hss\ \ 1\hss&\vtrule
\hss\ \ 1\hss&\vtrule
\hss\ \ 0\hss&\vtrule
\hss\ \ 1\hss&\vtrule
\hss\ \ 1\hss&\vtrule
\hss\ \ 1\hss&\vtrule
\hss\ \ 1\hss&\vtrule
\hss\ \ 0\hss&\vtrule
\hss\ \ 1\hss&\vtrule
\hss\ \ 1\hss&\vtrule
\cr\hrule
\+\vtrule
\hss\ p=9\hss&\vtrule
\hss\ \ 0\hss&\vtrule
\hss\ \ 1\hss&\vtrule
\hss\ \ 0\hss&\vtrule
\hss\ \ 1\hss&\vtrule
\hss\ \ 1\hss&\vtrule
\hss\ \ 1\hss&\vtrule
\hss\ \ 1\hss&\vtrule
\hss\ \ 1\hss&\vtrule
\hss\ \ 1\hss&\vtrule
\hss\ \ 0\hss&\vtrule
\hss\ \ 1\hss&\vtrule
\cr\hrule
\+\vtrule
\hss\ p=10\hss&\vtrule
\hss\ \ 0\hss&\vtrule
\hss\ \ 0\hss&\vtrule
\hss\ \ 1\hss&\vtrule
\hss\ \ 1\hss&\vtrule
\hss\ \ 1\hss&\vtrule
\hss\ \ 1\hss&\vtrule
\hss\ \ 1\hss&\vtrule
\hss\ \ 1\hss&\vtrule
\hss\ \ 1\hss&\vtrule
\hss\ \ 1\hss&\vtrule
\hss\ \ 0\hss&\vtrule
\cr\hrule
}}
\mytag{2.5}
$$
In the memory of a computer bit-arrays like \mythetag{2.5} are packed into
byte-arrays with $8$ bits per $1$ byte, e.\,g\. the array $u_{11}$ looks
like
$$
\hskip -2em
\vcenter{\vbox{\hsize 12cm
\offinterlineskip\settabs\+\indent
\vtrule
\hskip 1.2cm &\vtrule
\hskip 1.7cm &\vtrule
\hskip 3.8cm &\vtrule
\hskip 1.7cm &\vtrule
\cr\hrule
\+\vtrule
\hss\ $00000000$\hss&\vtrule
\hss\ $11100000$\ \hss&\vtrule
\hss\ $10011111\ \ldots\ $11111111$ $\hss&\vtrule
\hss\ $0\blue{0000000}$\ \hss&\vtrule
\cr\hrule
}}
\mytag{2.6}
$$
Note that bits in a byte are written in the reverse order --- the highest bit
is the leftmost. This is because bytes are designed to represent binary numbers. 
Note also that the last byte of the table \mythetag{2.5} in \mythetag{2.6} is
incomplete. It is appended with zero bits which are shown in blue.\par
     Bytes associated with the  prime number $r=11$ can be written into some
linear locus of memory. Similarly, bytes associated with several other prime 
numbers can be written into adjacent loci. Altogether they constitute a bit seive. 
Accessing a proper bit of this seive, we can easily decide whether for a certain 
pair of integer numbers $p$ and $q$ the equation $Q_{pq}(t)=0$ is unsolvable 
modulo some prime number $r$ enclosed in the seive. Then it is unsolvable in the 
ring of integers $\Bbb Z$ as well. Quickening the search algorithm is reached 
through sweeping away those $(p,q)$ pairs that do not go through the bit seive
for several prime numbers. Indeed, it is clear that calculating the remainders 
$$
\xalignat 2
&\tilde p=p\!\!\mod r,
&&\tilde q=q\!\!\mod r.
\endxalignat
$$
and then addressing bits in a memory locus are much faster operations than 
factoring a polynomial with numeric coefficients.\par
     In our particular case we use the bit seive for 96 consecutive prime numbers 
from $11$ to $541$. This bit seive is stored in the binary file 
\darkred{{\tt Cuboid\kern 0.5pt\_\kern 1pt pq\kern -0.5pt\_\kern 0.5pt bit\_\kern 
0.5pt tables.bin}}. In order to access effectively bit-tables for each particular 
prime number from $11$ to $541$ one should know their offsets within this file. 
These offsets are written to the separate binary file 
\darkred{{\tt Cuboid\kern 0.5pt\_\kern 0.5pt primes.bin}}. They are 
enclosed in the structures described as follows in C++ language:\par
\parshape 1 10pt 350pt 
{\tt\noindent
\blue{struct} primes\_\kern 1pt item\newline 
\{\newline 
\blue{short} prime;	\ \ \ \ \ \ \ \ \ \ \green{// prime number}\newline 
\blue{unsigned int} p\_\kern 1pt offset;	\green{// prime bit-table offset}\newline  
\};}\par
In our implementation the values of prime numbers are restricted not only by
{\tt short=2\,bytes} data format used for them. Each prime number $r$ is associated
with the $r\times r$ bit-table that occupies $r^{\kern 0.5pt 2}/8$ bytes in memory.
Using {\tt unsigned int=4\,bytes} format for offsets, we have the following restriction:
$$
\hskip -2em
\sum^N_{i=5}\frac{r^{\kern 0.5pt 2}_i}{8}< 2^{\kern 0.5pt 32}
\text{, \ where \ }r_5=11,\ r_6=13,\,\ldots\,.
\mytag{2.7}
$$
From \mythetag{2.7} we derive $N<1198$ and $r_N<9697$. These inequalities fit 
the 4\,Gb RAM (random access memory) limit. Actually we have chosen 
$N=100$ in which case 1.5\,Mb RAM is sufficient.\par
\head
3. Code for generating binary files. 
\endhead
     The code for preparing \darkred{{\tt Cuboid\kern 0.5pt\_\kern 
0.5pt pq\kern -0.5pt\_\kern 0.5pt bit\_\kern 0.5pt tables.bin}} and 
\darkred{{\tt Cuboid\kern 0.5pt\_\kern 0.5pt primes.bin}} binary files is implemented 
as a DLL library interacting with a Maple code. The DLL file 
\darkred{{\tt Cuboid\_\kern 0.5pt search\_\kern 0.5pt v01.dll}} is generated
within the 32-bit x86 makefile project for Microsoft Visual C++ 2005 Express 
Edition package. The project files
\roster
\item"1)"\quad\darkred{{\tt make.bat}}
\item"2)"\quad\darkred{{\tt makefile}}
\item"3)"\quad\darkred{{\tt Cuboid\_\kern 0.5pt search\_\kern 0.5pt v01.h}}
\item"4)"\quad\darkred{{\tt Cuboid\_\kern 0.5pt search\_\kern 0.5pt v01.cpp}}
\endroster
are suppled as ancillary files to this paper. 
The C++ file \darkred{{\tt Cuboid\_\kern 0.5pt search\_\kern 0.5pt v01.cpp}}
\linebreak is the main source file of the project. It comprises a C++ and 
inline assembly language code for running on 32-bit Windows machines with 
Intel compatible processors. The DLL library file \darkred{{\tt Cuboid\_\kern 
0.5pt search\_\kern 0.5pt v01.dll}} is produced from this code by running the 
batch file \darkred{{\tt make.bat}} in a command prompt window:
\medskip
$>$\,{\tt make.bat}
\medskip
\noindent Along with \darkred{{\tt Cuboid\_\kern 0.5pt search\_\kern 0.5pt 
v01.dll}} several other files are generated, including two LOG files 
\darkred{{\tt compiler.log}} and \darkred{{\tt linker.log}}. They can be removed
by running the same batch file with the {\tt clean} option:
\medskip
$>$\,{\tt make.bat clean}
\medskip
\noindent Note that the Visual C++ 2005 Express Edition package should be installed
for successfully running the above files. In our implementation it was installed
on Windows XP machine with Intel Pentium 4 Prescott CPU 2.80 GHz.\par
     The generated DLL library \darkred{{\tt Cuboid\_\kern 0.5pt search\_\kern 
0.5pt v01.dll}} exports several functions. Their declarations are in the C++ header 
file \darkred{{\tt Cuboid\_\kern 0.5pt search\_\kern 0.5pt v01.h}}. \pagebreak
Three of these functions are declared as follows:
\medskip
\parshape 1 10pt 350pt 
{\tt\noindent
\blue{extern} \darkred{"C"} 
\blue{\_\_declspec}(\blue{dllexport})\newline
\hphantom{a}\ \ \ \ \ \ \ \ \ \ \ \ \ \ \ \ \ \ \ \ \ \ \ \ 
\blue{void \_\_stdcall} Open\_\kern 0.5pt pq\_\kern 1pt file\_\kern 0.5pt stream();
\newline
\blue{extern} \darkred{"C"} 
\blue{\_\_declspec}(\blue{dllexport})\newline
\hphantom{a}\ \ 
\blue{unsigned int \_\_stdcall} Write\_\kern 0.5pt pq\_\kern 1pt file\_\kern 0.5pt 
stream(\blue{unsigned int} rrr);\newline
\blue{extern} \darkred{"C"} 
\blue{\_\_declspec}(\blue{dllexport})\newline
\hphantom{a}\ \ \ \ \ \ \ \ \ \ \ \ \ \ \ \ \ \ \ \ \ \ \ \ 
\blue{void \_\_stdcall} Close\_\kern 0.5pt pq\_\kern 1pt file\_\kern 0.5pt stream();
}
\medskip
\noindent
These declarations correspond to the following Maple worksheet declarations:
\medskip
\noindent{\tt\red{>\ DLL\_\kern 1pt file:="./Cuboid\_\kern 0.5pt search\_\kern 0.5pt
v01.dll":\newline
>\ Open\_\kern 0.5pt pq\_\kern 1pt file\_\kern 0.5pt stream:=define\_\kern 0.5pt
external('Open\_\kern 0.5pt pq\_\kern 1pt file\_\kern 0.5pt stream',\newline
\hphantom{a}\ \ \ LIB=DLL\_\kern 1pt file):\newline
>\ Write\_\kern 0.5pt pq\_\kern 1pt file\_\kern 0.5pt stream:=define\_\kern 0.5pt
external('Write\_pq\_\kern 1pt file\_\kern 0.5pt stream',\newline
\hphantom{a}\ \ \ 'rrr'::(integer[4]), RETURN::(integer[4]), LIB=DLL\_\kern 1pt 
file):\newline
>\ Close\_\kern 0.5pt pq\_\kern 1pt file\_\kern 0.5pt stream:=define\_\kern 0.5pt
external('Close\_\kern 0.5pt pq\_\kern 1pt file\_\kern 0.5pt stream',\newline
\hphantom{a}\ \ \ LIB=DLL\_\kern 1pt file):}}
\medskip
\noindent
Maple worksheets are supplied in XML format as ancillary files to this paper. 
They can be imported to Maple. Here is the list of these files: 
\roster
\item"5)"\quad\darkred{{\tt Create\_\kern 0.5pt binary\_\kern 0.5pt 
seive\_\kern 1pt files.xml}}
\item"6)"\quad\darkred{{\tt Test\_\kern 1pt external\_\kern 1pt 
DLL\_\kern 0.5pt procedures\_01.xml}}
\item"7)"\quad\darkred{{\tt Test\_\kern 1pt external\_\kern 1pt 
DLL\_\kern 0.5pt procedures\_02.xml}}
\item"8)"\quad\darkred{{\tt Search\kern 0.5pt\_\kern 1pt for\_\kern 0.5pt 
cuboids.xml}}
\endroster
\par 
     The external function \red{{\tt Write\_\kern 0.5pt pq\_\kern 1pt file\_\kern 
0.5pt stream(r)}} imported to the Maple worksheet creates the bit-seive table 
for a given prime number \red{r} in its argument and writes it to the binary 
file \darkred{{\tt Cuboid\kern 0.5pt\_\kern 0.5pt pq\kern -0.5pt\_\kern 0.5pt 
bit\_\kern 0.5pt tables.bin}}. It returns the integer value equal to the number of 
bytes written to the file \darkred{{\tt Cuboid\kern 0.5pt\_\kern 0.5pt pq\kern 
-0.5pt\_\kern 0.5pt bit\_\kern 0.5pt tables.bin}}. The other binary file 
\darkred{{\tt Cuboid\kern 0.5pt\_\kern 0.5pt primes.bin}} is written simultaneously 
using the Maple worksheet code in \darkred{{\tt Create\_\kern 0.5pt binary\_\kern 0.5pt 
seive\_\kern 1pt files.xml}}.\par
      The external function \red{{\tt Write\_\kern 0.5pt pq\_\kern 1pt file\_\kern 
0.5pt stream(r)}} exploits another external function 
\red{{\tt Calculate\_Q\_\kern 0.5pt pq\_\kern 0.5pt mod\_\kern 0.5pt prime(p,q,t,r)}}. 
This function returns the value of the polynomial \mythetag{1.1} modulo prime 
number \red{r} taken as its fourth argument. Though its arguments are declared as 
32-bit integers, its code is designed to deal with 16-bit unsigned integers only. 
Due to the restriction $r<9697$ derived from \mythetag{2.7} its usage in 
\red{{\tt Write\_\kern 0.5pt pq\_\kern 1pt file\_\kern 0.5pt stream(r)}} does not 
require 32-bit integers in its arguments: 
$$
\xalignat 2
&0\leqslant p\leqslant r-1<9797,
&&0\leqslant q\leqslant r-1<9797.
\endxalignat
$$
The function \red{{\tt Calculate\_Q\_\kern 0.5pt pq\_\kern 0.5pt mod\_\kern 0.5pt
prime(p,q,t,r)}} is a delicate part of the pro\-ject. It is written in assembly language. 
Therefore it is carefully tested in the Maple worksheet code file 
\darkred{{\tt Test\_\kern 1pt external\_\kern 1pt DLL\_\kern 0.5pt 
procedures\_01.xml}}.\par
\head
4. Code for loading and unloading binary files. 
\endhead
    Once the binary files \darkred{{\tt Cuboid\kern 0.5pt\_\kern 
0.5pt pq\kern -0.5pt\_\kern 0.5pt bit\_\kern 0.5pt tables.bin}} and 
\darkred{{\tt Cuboid\kern 0.5pt\_\kern 0.5pt primes.bin}} are generated, they should 
be used in searching for perfect cuboids. For this purpose they should be loaded into 
the memory easily accessible from the DLL library functions. This task is performed 
by the function \red{{\tt Load\_Cuboid\_\kern 1pt Binaries()}} residing within the 
same DLL library. The opposite task is to unload the binary files, i\.\,e\. to 
release the memory occupied by them. This task is performed by the function
\red{{\tt Release\_Cuboid\_\kern 1pt Binaries()}} also residing within the DLL 
library.\par 
     The loading and unloading functions are imported and tested within the 
Maple worksheet code file \darkred{{\tt Test\_\kern 1pt external\_\kern 1pt 
DLL\_\kern 0.5pt procedures\_02.xml}}.\par
     As an auxiliary test for two generated binary files \darkred{{\tt Cuboid\kern 
0.5pt\_\kern 0.5pt pq\kern -0.5pt\_\kern 0.5pt bit\_\kern 0.5pt tables.bin}} and 
\darkred{{\tt Cuboid\kern 0.5pt\_\kern 0.5pt primes.bin}} we visualize the 
bit-tables from \darkred{{\tt Cuboid\kern 0.5pt\_\kern 0.5pt pq\kern 
-0.5pt\_\kern 0.5pt bitmaps.bin}} in the form of the text file \darkred{{\tt 
Cuboid\kern 0.5pt\_bit\_\kern -0.5pt tables.txt}}. This text file is written by 
the code from the Maple worksheet file \darkred{{\tt Test\_\kern 1pt external\_\kern 
1pt DLL\_\kern 0.5pt procedures\_02.xml}}.\par
\head
5. Code for searching cuboids. 
\endhead
     According to Theorem~\mythetheorem{1.2} the search for cuboids in the case of
the second cuboid conjecture consists in scanning the region given by the
inequalities \mythetag{1.2}. For each positive $p$\/ these inequalities specify a finite
segment of the real axis, which comprises a finite number of integer points. For 
$p\leqslant 151$ this segment is given by
$$
\hskip -2em
\frac{p}{59}\leqslant q\leqslant{59\,p}.
\mytag{5.1}
$$
For $p\geqslant 152$ the inequalities are different:
$$
\hskip -2em
\root{\raise 1pt\hbox{$\ssize 3\kern -1pt $}}\of{\frac{p}{9}\,}
\leqslant q\leqslant{59\,p}.
\mytag{5.2}
$$
The inequalities $0<p\leqslant 151$ and the inequalities \mythetag{5.1} 
outline a finite set of integer points on the coordinate $pq$-plane. One 
can easily verify that these points do not produce perfect cuboids. For 
this reason the software in the DLL library \darkred{{\tt Cuboid\_\kern 
0.5pt search\_\kern 0.5pt v01.dll}} is designed to search cuboids 
for $p\geqslant 152$ within each segment specified by the inequalities 
\mythetag{5.2}. Roughly speaking, it is an infinite loop on $p\geqslant
152$ and an enclosed loop on $q$ obeying the inequalities \mythetag{5.2} 
for each $p$.\par
     Both loops on $p$ and on $q$ are started from within the Maple
worksheet file \darkred{{\tt Search\kern 0.5pt\_\kern 1pt for\_\kern 0.5pt 
cuboids.txt}} by executing the commands 
\medskip
\noindent{\tt\red{>\ Load\_Cuboid\_\kern 1pt Binaries();\newline
>\ Start\_\kern 1pt searching(152,3);
}}
\medskip
\noindent
Here $152$ and $3$ are initial values for $p$ and $q$ respectively. They 
should obey the inequalities \mythetag{5.2}. The function 
\red{{\tt Start\_\kern 0.5pt searching}} is an external function imported 
from the DLL library \darkred{{\tt Cuboid\_\kern 0.5pt search\_\kern 0.5pt 
v01.dll}}. It starts the looping process and returns just immediately with 
the value \red{{\tt 0}} indicating that the search is successfully started.
The multithreading mechanism is used in the code of this function:
\medskip
\parshape 1 10pt 350pt
\noindent 
{\tt\_\kern 0.5pt beginthread(Look\_\kern 0.5pt for\_\kern 0.5pt
cuboids\_thread,0,(\blue{void}*)12);\newline
\blue{return}(0); 
}
\medskip
\noindent 
Here \darkred{{\tt Look\_\kern 0.5pt for\_\kern 0.5pt
cuboids\_thread}} is an internal function which is not exported from
the DLL library. It is executed within a new thread, while the initial function
\darkred{{\tt Start\_\kern 0.5pt searching}} returns control to the Maple 
worksheet.\par
     You can do anything in the Maple worksheet while the search function 
\darkred{{\tt Look\_\kern 0.5pt for\_\linebreak \kern 0.5pt cuboids\_thread}} 
is running its infinite loops on $p$ and $q$, provided you do not stop the Maple 
session by closing the worksheet. \pagebreak In particular, you can control the 
process of searching by executing the following function in the Maple worksheet:
\medskip
\noindent{\tt\red{>\ Get\_\kern 0.5pt current\_\kern 1pt p();
}}
\medskip
\noindent
This function returns the current value of the loop variable $p$ with the 
discreteness equal to $100$. There is another control function: 
\medskip
\noindent{\tt\red{>\ Get\_\kern 0.5pt current\_\kern 1.5pt r\kern -0.5pt 
\_\kern 1.5pt max();
}}
\medskip
\noindent
This function returns the maximal prime number is used to seive cuboids within
current hundred values of $p$. The number $r_{\text{max}}$ is flushed to $1$
for each next hundred values of $p$ and then is recalculated again.\par
    The infinite loop on $p$ cannot ever terminate by itself. Therefore it is
terminated manually. This can be done at any time by executing the following 
function in the Maple worksheet that initiated the thread with this loop:
\medskip
\noindent{\tt\red{>\ Stop\_\kern 1pt searching();
}}
\medskip
\noindent
Upon doing it you can read the time stamp and the exit values of $p$, $q$, and
$r_{\text{max}}$ at the end of the file \darkred{{\tt Cuboid\_\kern 0.5pt 
search\_\kern 1pt report.txt}}, e\.\,g\. it could be
\medskip
\noindent{\tt
\green{2016-2-20 21:36\newline 
Stop with p=112618, q=5691455, r\kern -0.5pt \_\kern 1.5pt max=131
}}
\medskip
\noindent
Then you can restart the search from this point on by executing the command
\medskip
\noindent{\tt\red{>\ Start\_\kern 1pt searching(112618,5691455);
}}
\medskip
\noindent
Or you can terminate the session by executing the command
\medskip
\noindent{\tt\red{>\ Release\_\kern 0.5pt Cuboid\_\kern 1pt Binaries();
}}
\medskip
\noindent
and then closing the Maple worksheet.\par
    Normally the function  \red{{\tt Start\_\kern 1pt searching(p,q)}} returns
$0$ indicating that the search is successfully started. It returns $1$ if the
search is already running. So you cannot initiate several search threads running
simultaneously with this software. This limitation is planned to be removed in 
further versions of the DLL library \darkred{{\tt Cuboid\_\kern 0.5pt search\_\kern 
0.5pt v01.dll}}.\par
    The function \red{{\tt Start\_\kern 1pt searching(p,q)}} returns $2$ if 
$p<152$ (see \mythetag{5.1} and \mythetag{5.2} above for explanation). This 
function returns $3$ if $p>72796055$, which breaks the 32-bit limit for $q=59\,p$ 
in \mythetag{5.2}.\par
    The function \red{{\tt Start\_\kern 1pt searching(p,q)}} returns $4$ if
$q$ is below the lower limit set by the inequalities \mythetag{5.2}. Similarly,
it returns $5$ if $q$ is above the upper limit set by the inequalities 
\mythetag{5.2}.\par
    The function \red{{\tt Start\_\kern 1pt searching(p,q)}} returns $6$ if it
is invoked before the binary files \darkred{{\tt Cuboid\kern 
0.5pt\_\kern 0.5pt pq\kern -0.5pt\_\kern 0.5pt bit\_\kern 0.5pt tables.bin}} 
and \darkred{{\tt Cuboid\kern 0.5pt\_\kern 0.5pt primes.bin}} are loaded into
the memory by the function \red{{\tt Load\_Cuboid\_\kern 1pt Binaries()}}.
\par
    The function \red{{\tt Load\_Cuboid\_\kern 1pt Binaries()}} normally returns
the number of bytes loaded from the file \darkred{{\tt Cuboid\kern 0.5pt\_\kern 
0.5pt primes.bin}}, i\.\,e\. the size of this file. However, if it is invoked
when the binary files \darkred{{\tt Cuboid\kern 
0.5pt\_\kern 0.5pt pq\kern -0.5pt\_\kern 0.5pt bit\_\kern 0.5pt tables.bin}} and 
\darkred{{\tt Cuboid\kern 0.5pt\_\kern 0.5pt primes.bin}} are already loaded, it 
does not load them again and returns $0$.\par 
    The function \red{{\tt Release\_\kern 0.5pt Cuboid\_\kern 1pt Binaries()}}
normally returns $0$. However it returns $6$ if the binary files \darkred{{\tt 
Cuboid\kern 0.5pt\_\kern 0.5pt pq\kern -0.5pt\_\kern 0.5pt bit\_\kern 0.5pt 
tables.bin}} and \darkred{{\tt Cuboid\kern 0.5pt\_\kern 0.5pt primes.bin}} are
not loaded. This function cannot release the memory occupied by the bit-tables 
if the search is running. In this case it returns $1$.\par 
    The function \red{{\tt Stop\_\kern 1pt serching()}} normally returns $0$.
However, if this function is invoked when the search is not on, it returns $1$. 
Thus the functions
\medskip
\noindent{\tt\red{>\ Load\_Cuboid\_\kern 1pt Binaries();\newline
>\ Start\_\kern 1pt searching(p,q);\newline
>\ Stop\_\kern 1pt serching();\newline
>\ Release\_\kern 0.5pt Cuboid\_\kern 1pt Binaries();
}}
\medskip
\noindent    
should be invoked in the above order. Otherwise they signal misuse, but do
not lead to a crash. These four functions constitute a toolkit we used in the 
present numerical research of perfect cuboids.\par  
\head 
6. Results. 
\endhead 
     At present date 01.04.2016 the values of $p$ from 1 to 154000 are scanned. For 
each such $p$ all values of $q$ limited by the inequalities \mythetag{1.2} are scanned. 
This stands for about $700$ billions $(p,q)$ pairs that have been tested. Indeed, we 
have 
$$
N\approx\!\!\sum^{154000}_{p=1}\!\!\!59\,p=699626543000\approx 0.7\cdot 10^{12}.
$$
None of these $(p,q)$ pairs produces a perfect cuboid. Moreover, none of them goes
through our primes seive composed by 96 consecutive primes from 11 to 541. Actually 
this seive is so dense that the maximal depth reached thus far is 29th prime in our 
sequence, which is equal to 137. This result is negative in the sense of finding
a perfect cuboid. However it shows that the Second cuboid 
conjecture~\mytheconjecture{1.1} is rather firm for to believe that it might be 
valid.\par
     The total time spent for the above computations is 4130 minutes, i\.\,e\. 
about 69 hours. Then we can calculate the time per one $(p,q)$ pair:
$$
\Delta\,t=\frac{4130}{699626543000}\,\text{\sl min}\approx 3.54\cdot 10^{-7}\,
\text{\sl sec}.
$$
The upper limit for $p$ with our 32-bit code is given by the formula
$$
p_{\,\text{max}}=\frac{2^{\,32}}{59}\approx 72796055. 
$$
Here is the estimate for the time needed to reach this limit:
$$
t=\Delta\,t\kern -0.5em\sum^{72796055}_{p=1}\!\!\!59\,p=5.53\cdot
10^{10}\,\text{\sl sec}\approx 1755\,\text{\sl  years}. 
$$
This estimate means that our code should be improved not only at the expense 
of multithreading and multiprocessing. Some fresh theoretical ideas are 
required. 
\head
7. Acknowledgement.
\endhead
     The authors are grateful to A\.\,A\.\,Gubarev for helpful advices on
linking C++ code with a Maple worksheet session. 

\enddocument
\end